\documentclass[12pt]{article}
\font \bb=msbm10 %scaled \magstep2 % Blackboard bold for real numbers field
\def\C{\hbox{\bb C}}
\def\K{\hbox{\bb K}}
\begin{document}
\begin{flushleft}
\Large { \bf Generalized Symmetries
of Partial Differential Equations and Quasiexact Solvability}
\newtheorem{definition}{Definition}
\newtheorem{Proposition}{Proposition}
\newtheorem{Lemma}{Lemma}
\newtheorem{Theorem}{Theorem}
\end{flushleft}
\begin{flushleft}
{\it Arthur SERGHEYEV} \footnote{e-mail: arthur@apmat.freenet.kiev.ua}
\end{flushleft}
\begin{flushleft}
{\it Institute of Mathematics of the National Academy of Sciences of
Ukraine, \\
3 Tereshchenkivs'ka Street, Kyiv 4, \ Ukraine} \\
\end{flushleft}
\begin{minipage}{12cm}
\footnotesize

Using the adjoint action of the infinitesimal translations (with
respect to some (in)dependant variables) on specific
finite-dimensional subspaces of the space of generalized symmetries
of some system of partial differential equations, we explicitly
determine the dependance of coefficients of generalized symmetries
from these subspaces on the above-mentioned variables.
We establish the connection of our results with the theory of
quasiexactly solvable models.
Some generalizations of the approach proposed also are discussed.
\end{minipage}
%\sloppy
\section*{Introduction}

It is well known \cite{o}, \cite{fn} that the problem of finding of all the
generalized symmetries (or non-Lie symmetries in terms of \cite{fn})
of a given system of partial differential equations (PDEs)
is nontrivial and seldom admits the complete solution. More or less
complete results in this field are obtained for some linear equations
(e.g. free Dirac, Klein -- Gordon and Schr\"odinger equations (\cite{fn} and
references therein)),
1+1-dimensional evolution equations (vide, e.g., \cite{i}, \cite{s}) and
integrable systems (\cite{kp} and references therein).

On the other hand, the knowledge of the generalized symmetries of a
given system of PDEs allows one to construct (partially) invariant solutions
of this system (the famous examples of which are finite-gap solutions
of the KdV equation and numerous spherically-symmetric, self-similar
etc. solutions \cite{i}), to find the conservation laws \cite{o},
\cite{i}, to check the integrability \cite{s}
and to separate variables in linear PDEs \cite{m}. Thus, finding
generalized symmetries of given system of PDEs is interesting problem of
modern mathematical physics.

In this paper we consider the systems of PDEs, whose set of generalized
symmetries possesses specific finite-dimensional subspaces $V^{(q)}$
(refer to Section 2), which
are invariant under the action of the set of
infinitesimal translations with respect to some (in)dependant variables.
For the generalized symmetries from these subspaces we determine their
dependance on these variables, i.e. partially solve the
above-mentioned problem.

 The plan of the paper is as follows. In Section 1 we briefly remind
some basic notions concerning generalized symmetries. Section 2 presents
the key idea of the proposed approach and some possible
generalizations. Finally, in the Section 3 we discuss the connection
of our results with the theory of quasiexact solvability.

\section{Generalized symmetries: basic definitions}
Let us consider the system of PDEs of the form:
\begin{equation}  \label{1}
F _{\nu}(x,u,\dots, u^{(d)}) = 0, \quad \nu =1, \dots ,f,
\end{equation}
where $u=u(x)$ is unknown vector-function $u = (u_1, \dots, u _n)^{T}$ of $m$
independant variables $x=(x_1, \dots, x_m)$;
%$F = (F_1,\dots, F_f)$;
$u^{(s)}$
denotes the set of derivatives of $u$ with respect to $x$ of the order $s$.
( $^{T}$ denotes here and below matrix transposition).

\begin{definition}$\lbrack {\rm 1} \rbrack$. The operator ${\rm Q}$
of the form
\begin{equation} \label{2}
{\rm Q} = \sum _{i=1} ^{m} \xi _i (x,u,\dots, u^{(q)}) \partial
/\partial x_i +
\sum _{\alpha =1} ^{n} \eta _{\alpha} (x,u,\dots, u^{(q)}) \partial /
\partial u_{\alpha}
\end{equation}
is called the generalized symmetry of order $q$ of the system of PDEs
(\ref{1}) if
its prolongation ${\rm {\bf pr}}{\rm Q}$ annulates (\ref{1}) on the set $M$ of
(sufficiently smooth) solutions of (\ref{1}):
\begin{equation} \label{3}
  {\rm {\bf pr} Q} [F_ {\nu}] \mid _M = 0, \quad  \nu =1, \dots ,f.
\end{equation}
\end{definition}
Here \cite{o}
\begin{equation}  \label{4}
  {\rm {\bf pr} Q} = \sum _{l=1} ^{m} \xi _l D_l +
\sum _{\alpha =1} ^{n} \sum _{J}
%\in {\rm {\bf Z}}^{m}_{+} }
 (D _J (\eta _{\alpha}
- \sum _{l=1}^{m} \xi _l u_{\alpha ,l}) ) \frac{\partial}{\partial
u_{\alpha , J} } ,
\end{equation}
where the summation is extended to all the multiindices $J = (j_1,
\dots ,j _m)$ with non-negative integer $j_s$, $s=1, \dots , m$,
% \in {\rm {\bf Z}}_{+}$
%(${\rm {\bf Z}}_{+}$ is the set of non-negative integers),
$u_{\alpha , J} =\partial ^{j_1 +\dots
+j_m} u_{\alpha} / \partial x_1^{j_1} \partial x_{2} ^{j_2} \dots
\partial x_{m} ^{j _m}$, $u_{\alpha ,l} \equiv \partial u_{\alpha} /
\partial x_l$, $D_l$ is the so-called total derivative \cite{o} acting on
the functions of $x, u, u ^{(1)}, \dots $ :
\begin{displaymath}
D_l = \partial /\partial x_{l} + \sum _{J}
%\in {\rm {\bf Z}}_{+}^{m} }
(\partial u_{\alpha ,J}/\partial x_{l}) \frac{\partial}{\partial
u_{\alpha ,J}}
\end{displaymath}
and $D_J = D_1 ^{j_1} D_2 ^{j_2} \dots D_m ^{j_m}$.

Let us mention that the generalized symmetries of (\ref{1}) of order 0 are
usually called Lie symmetries \cite{o}.

\begin{definition} {\rm \cite{o} .}
Let ${\rm Q}_s$ be some differential operators of the form (\ref{2}):
\begin{displaymath}
{\rm Q}_s = \sum _{i=1} ^{m} \xi _{i} ^{(s)} (x,u,\dots, u^{(q_s)}) \partial /
\partial x_i +
\sum _{\alpha =1} ^{n} \eta _{\alpha} ^{(s)} (x,u,\dots, u^{(q_s)})
\partial / \partial u_{\alpha}
\end{displaymath}
Then the operator ${\rm Q}_3$ of the same form with the coefficients
\begin{displaymath}
\begin{array}{lll}
\xi _{i} ^{(3)} =
{\rm {\bf pr} Q}_1 [\xi _{i} ^{(2)}] - {\rm {\bf pr} Q}_2 [\xi _{i}
^{(1)}] \\
\eta _{\alpha} ^{(3)} = {\rm {\bf pr} Q}_1 [\eta _{\alpha} ^{(2)}]
- {\rm {\bf pr} Q}_2 [\eta _{\alpha} ^{(1)}]
\end{array}
\end{displaymath}
is called the Lie bracket of the operators ${\rm Q} _1$, ${\rm Q}_2$
and is denoted ${\rm Q}_{3} = [{\rm Q} _{1}, {\rm Q}_{2}]$.
\end{definition}
\begin{Proposition} {\rm \cite{o}}. If ${\rm Q}_1, {\rm Q}_2$ are generalized
symmetries of (\ref{1}), then so does $[{\rm Q}_1, {\rm Q}_2]$.
\end{Proposition}

 Let $Sym$ be the Lie algebra \cite{o}(with respect to Lie bracket $[,]$)
of all the generalized symmetries of (\ref{1}) of non-negative orders,
$Sym ^{(q)}$ the linear space of the generalized symmetries of (\ref{1})
of order not higher than $q$ and $p^{(q)}$
its dimension, $Sym _{q} \equiv Sym^{(q)}/Sym^{(q-1)}$ $(q \neq 0)$,
$Sym _0 \equiv Sym^{(0)}$ and $p _{q}$ the dimension of $Sym_q$.
It is is straightforward to check \cite{o} that $Sym _0$ is
subalgebra of $Sym$
with respect to Lie bracket. The generalized symmetries of (\ref{1}) from
$Sym _0$ (i.e. Lie symmetries) may be considered as vector fields on
the manifold of 0-jets $M^{(0)}$ \cite{o} with local coordinates
$z_A$: $z_{i} = x_{i}, i=1, \dots, m, z_{m+\alpha } = u _{\alpha},
\alpha =1, \dots , n$ (from now on the indices $A, B, C, D, \dots$
will run from
1 to $m+n$ and we shall denote $\partial _A \equiv \partial /\partial z_A$).
Moreover, the Lie bracket of two Lie symmetries of (\ref{1}) ${\rm Q}_{1},
{\rm Q}_{2}$ coincides with the commutator of corresponding vector fields.

Let us mention that we
choose as the basic field the field of complex numbers $\C$, i.e.
$Sym$ is considered as Lie algebra over
$\C$ and the coefficients $\xi _i, \eta _{\alpha}$ of generalized
symmetries are complex \footnote{In fact, everywhere in our considerations
(except the examples) $\C$ may be replaced by the arbitrary
algebraically closed field $\K$, since the Corollary 1 of the Theorem
$3'$ from $\S 2$ of Chapter VIII of \cite{g}, which we use in the
proofs of our theorems, remains true
for this case too .}.
\section{The explicit formulas for the symmetries}

\begin{Theorem} \label{t4}
Let $W$ be some linear subspace of the linear space of the
differential operators of the form (\ref{2}) of all orders
\footnote{$A_1, \dots , A_g$ are fixed integers from the range $1,
\dots, m+n$, $g \leq m+n$.},
$W_{A_{1}, \dots, A_{g}}$ be the linear space of the operators, obtained from
the operators from $W$ by setting $z_{A_{1}}=0, \dots, z_{A_g} = 0$ in their
coefficients, $V= W \bigcap Sym$, for some $q_1$ the dimension of the
subspace of $V$ $V^{(q_{\lower2pt\hbox{\tiny{1}}})} \equiv V \bigcap Sym
^{(q_{\lower2pt\hbox{\tiny{1}}})}$
$v^{(q_{\lower2pt\hbox{\tiny{1}}})} < \infty$
and for any generalized symmetry of (\ref{1}) ${\rm Q} \in
V^{(q_{\lower2pt\hbox{\tiny{1}}})}$ $\partial {\rm Q} / \partial z_{A_s}
 \in V ,s=1, \dots, g$ \footnote{it is clear that this implies
$\partial {\rm Q} / \partial z_{A_s} \in V^{(q_{\lower2pt\hbox{\tiny{1}}})},
s=1, \dots, g$.}.

Then in each $V^{(q)}$ ($q=0,
\dots, q_1)$ there exists such a basis of linearly independant
generalized symmetries ${\rm Q}_l^{(q,\gamma)}$,
$l=1, \dots, r_{\gamma}^{(q)}$, $\gamma=1,\dots, \rho ^{(q)}$ ($\rho^{(q)}
\leq v^{(q)}$, $\sum _{\gamma =1}^{\rho^{(q)}} r_{\gamma}^{(q)} = v^{(q)}$)
that the generalized symmetries from this basis are given by the
formulas
\begin{equation} \label{8aaa}
\begin{array}{lll}
{\rm Q}_l^{(q,\gamma)} =
\exp(\sum _{s=1}^{g}
\lambda _{\gamma} ^{(q,A_{s})} z_{A_{s}}) \times \\
\times \sum _{j_{1}=0} ^{k_{\gamma}^{(q,A_{1})} - 1}
\dots \sum _{j_{g}=0} ^{k_{\gamma}^{(q,A_{g})} - 1}
%\prod _{h=1}^{g}
(z_{A_{1}}) ^{j_{1}} (z_{A_{2}}) ^{j_{2}} \dots (z_{A_{g}}) ^{j_{g}}
{\rm C}_{l,j_{1}, \dots , j_{g} }^{(q,\gamma)}
\end{array}
\end{equation}
where ${\rm C}_{l,j_{1}, \dots , j_{g} }^{(q,\gamma)}$
are some differential operators from $W_{A_{1},\dots, A_{g}}$
of order $q$ or lower;
%%with the coefficients independant from $z_{A_{1}}, \dots, z_{A_{g}}$;
$\lambda _{\gamma} ^{(q,A_{s})} \in \C$ are some constants and
$k_{\gamma}^{(q,A_{s})}, s=1, \dots, g$ are some fixed numbers from the
range $1,\dots, r_{\gamma}^{(q)}$.
\end{Theorem}
{\it Proof.} Let ${\rm Q}_{1} ^{(s)}, \dots,{\rm Q} _{v_s} ^{(s)}$ be
some basis in $V _s$, where $V_s= V^{(s)}/V^{(s-1)}, s \neq 0, V_{0}
=V^{(0)}$.

According to the conditions of the Theorem $\partial {\rm
Q}_{l}^{(s)}/\partial
z_{A_i}$ also are generalized symmetries of (\ref{1}), which
obviously belong to $V^{(s)}$ (but not necessarily to $V_s$).
Thus, the set of differential
operators $\partial /\partial z_{A_{i}}$ possesses
finite-dimensional
\footnote{since $v^{(s)} \leq v^{(q_{1})}$ for $s \leq q_{1}$.}
invariant spaces $V^{(s)}$, $s=0, \dots, q_{1}$. Let us denote the
finite-dimensional linear operator being the representation
of $\partial /\partial z_{A_{i}}$ on $V^{(s)}$ by ${\rm G}^{(s,
A_{i})}$. Since $\partial ^{2} {\rm Q}_{l}^{(s)} /\partial z_{A_{i}} \partial
z_{A_{j}}  =\partial ^{2} {\rm Q}_{l}^{(s)} /\partial z_{A_{j}} \partial
z_{A_{i}}$, i.e. the operators $\partial /\partial z_{A_{i}}$,
 $\partial /\partial z_{A_{j}}$ commute, so do their representations:
\begin{equation} \label{comm}
{\rm G}^{(s,A_{i})}{\rm G}^{(s,A_{j})} = {\rm G}^{(s,A_{j})}{\rm
G}^{(s,A_{i})},\quad i,j=1,\dots,g, s=0, \dots, q_{1}.
\end{equation}

 Hence, according to the Corollary 1 of the Theorem $3'$ from
$\S 2$ of Chapter VIII of \cite{g}, the space $V ^{(q)}$ ($q \leq q_{1}$)
may be decomposed into the direct sum of such common invariant spaces
$I_{\gamma}^{(q)}$ of the linear operators ${\rm G}^{(q, A_1)},
\dots, {\rm G}^{(q, A_g)}$, that the minimal polynomials of ${\rm
G}^{(q, A_s)}$ on $I_{\gamma}^{(q)}$ are
\begin{equation} \label{chareq}
({\rm G}^{(q,A_s)} - \lambda_{\gamma}^{(q,A_s)})^{k_{\gamma}^{(q,A_s)}} =0
\end{equation}
for some $\lambda_{\gamma}^{(q,A_s)} \in \C$ and
$k_{\gamma}^{(q,A_s)}$ ($k_{\gamma}^{(q,A_s)}$ are fixed integers
from the range $1, \dots, r_{\gamma}^{(q)}$, where $
r_{\gamma}^{(q)}$ denotes the dimension of the space  $I_{\gamma}^{(q)}$),
$\gamma=1, \dots, \rho^{(q)}$, $s=1, \dots, g$. In order to avoid possible
ambiguity in the
definition of the subspaces $I_{\gamma}^{(q)}$ we set the following
requirement: the representation  ${\rm G}_{\gamma}^{(q,A_i)}$
of the operators ${\rm G}^{(q,A_i)}$ on each $I_{\gamma}^{(q)}$,
$i=1, \dots, g$, must be indecomposable.

Using this result allows us to resrict ourselves to considering the
only $r_{\gamma}^{(q)}$-dimensional subspace $I_{\gamma}^{(q)}$ of
$V^{(q)}$ and some basis ${\rm Q}_l^{(q,\gamma)}$, $l=1, \dots,
r_{\gamma}^{(q)}$ in it.
Let ${\rm R}^{(q,\gamma)} \equiv ({\rm Q}_{1}^{(q,\gamma)}, {\rm
Q}_{2}^{(q,\gamma)}, \dots, {\rm Q}_{r_{\gamma}^{(q)}}^{(q,
\gamma)})^{T}$ and let ${\rm
G}_{\gamma}^{(q,A_s)}$ denote the restriction of ${\rm G}^{(q,A_s)}$ on
$I_{\gamma}^{(q)}$.

Then we obtain
\begin{displaymath}
\partial {\rm R}^{(q, \gamma)} /\partial z_{A_i} = {\rm
G}_{\gamma}^{(q,A_{i})} {\rm R} ^{(q, \gamma)},
\end{displaymath}
where we have identified the operator ${\rm G}_{\gamma}^{(q,A_i)}$
with its matrix in the basis ${\rm Q}_l^{(q,\gamma)}$,  $l=1, \dots,
r_{\gamma}^{(q)}$;
each such system is compatible, since the matrices ${\rm
G}_{\gamma}^{(q,A_i)}$ commute in virtue of (\ref{comm}), and its
general solution is:
\begin{equation}  \label{10b}
{\rm R}^{(q,\gamma)} =  \prod _{i=1}^{g} \exp({\rm G}_{\gamma}
^{(q,A_i)} z_{A_i})  {\rm C} ^{(q,\gamma)},
\end{equation}
where ${\rm C}^{(q,\gamma)}$ is the vector of the operators from $W$
of order $q$ or lower, whose coefficients are
independant from $z_{A_1},\dots, z_{A_g}$ (i.e.
${\rm C}^{(q,\gamma)} \in W_{A_{1}, \dots,A_{g}}$). In virtue of
(\ref{chareq}) we obtain
\begin{equation} \label{expg}
\begin{array}{lll}
\exp({\rm G}_{\gamma}^{(q,A_i)} z_{A_i}) \equiv exp( \lambda_{\gamma}^{(q,
A_i)} z_{A_i}) \exp(({\rm G}_{\gamma}^{(q,A_i)} - \lambda_{\gamma}^{(q,
A_i)}) z_{A_i}) = \\
=\exp( \lambda_{\gamma}^{(q,A_i)} z_{A_i})
\sum_{s=0}^{k_{\gamma}^{(q,A_i)}-1} \frac{(z_{A_i})^{s}}{s!} ({\rm
G}_{\gamma}^{(q,A_s)} - \lambda_{\gamma}^{(q,A_i)})^{s}.
\end{array}
\end{equation}
The substitution of (\ref{expg}) into (\ref{10b}) evidently yields
(\ref{8aaa}). $\triangleright$

{\it Remark 1.} Often (e.g. if $V$ is an ideal in $Sym$ or if $V$ is
Lie subalgebra of $Sym$, containing $\partial /\partial z_{A_{i}}$ or
generalized symmetries, which are equivalent to them)
the sufficient condition
for $\partial {\rm Q} / \partial z_{A_s} \in V, s=1, \dots, g$
%for $A_1, \dots, A_g \in \lbrace 1, \dots ,m+n\rbrace$
(if ${\rm Q} \in
V^{(q_{\lower2pt\hbox{\tiny{1}}})}$) is (in virtue of the Proposition
1) the existence of Lie
symmetries $\partial /\partial z_{A_1}, \partial
/\partial z_{A_2}, \dots ,\partial /\partial z _{A_g}$
for the system (\ref{1}).

{\it Remark 2.} In the above Theorem it is
possible to set $q_1 = \infty$, taking into account that by construction
$V ^{(\infty)} \equiv V$ and modifying its conditions in the following way:
$v^{(q_1)} < \infty$ is replaced by $v_{q} < \infty$ for all
$q=0,1,\dots$, where $v_{q}$ is the dimension of $V_q$.

{\it Remark 3.} If $W$ is the space of all the Lie vector fields on
$M^{(0)}$, $V \equiv W \bigcap Sym = Sym ^{(0)}$, $v^{(0)} < \infty$
and the system (\ref{1}) admits Lie symmetries $\partial /\partial
z_{A}, A=1, \dots, m+n$ the Theorem \ref{t4} allows us to find the
dependance of {\it all} the Lie symmetries on {\it all} the variables $x,u$
they depend from, i.e. to reduce the problem of description of
Lie symmetries of (\ref{1}) to solving algebraic equations
\footnote{Let us mention that in this point our results are
completely analogous to those of the theory of quasiexact solvability
(see also Section 3).}.

From the merely technical point of view our results mean that if the
conditions of
the Theorem \ref{t4} are valid, one may seek for {\it all} the generalized
symmetries of (\ref{1}) from $V^{(q)}$ ($q \leq q_{1}$) without loss
of generality in the form
\begin{equation} \label{ansatzzz}
\begin{array}{lll}
{\rm Q} = \exp(\sum _{s=1}^{g}
\lambda ^{(q,A_{s})} z_{A_{s}}) \times \\
\times \sum _{j_{1}=0} ^{v^{(q)} - 1}
\dots \sum _{j_{g}=0} ^{v^{(q)} - 1}
(z_{A_{1}}) ^{j_{1}} (z_{A_{2}}) ^{j_{2}} \dots (z_{A_{g}}) ^{j_{g}}
{\rm C}_{j_{1}, \dots , j_{g} },
\end{array}
\end{equation}
where $\lambda^{(q,A_s)} \in \C$ and ${\rm C}_{j_{1}, \dots , j_{g}}$ are
differential operators from $W_{A_{1}, \dots, A_{g}}$ of order $q$ or lower.
Thus, the substitution of (\ref{ansatzzz}) into (\ref{3}) yields the
equations for the coefficients of ${\rm C}_{j_{1}, \dots , j_{g}}$,
and if one is able to find all the independant solutions of these
equations, the substitution of them into (\ref{ansatzzz}) gives {\it
all} the linearly independant symmetries of (\ref{1}) from $V^{(q)}$.

Let us illustrate these ideas by the following examples:

{\it Example 1.} Let $m=2$, $n=1$, $u_1 \equiv u$, $x=(x_1 \equiv t,
x_2 \equiv y)$, $u _{(l)} = \partial ^{l} u / \partial y ^{l}$ and
$u$ satisfies the evolution equation
\begin{equation} \label{eveq}
\partial u / \partial t =G(u, u_{(1)}, \dots, u _{(d)}), \quad d\geq 2,
\end{equation}
$W$ be the linear space of the differential operators of the form (\ref{2})
with $\xi _i \equiv 0$, whose coefficient
$\eta \equiv \eta _{1}$, which is called the characteristics of the symmetry
\cite{s},
depends only on $y, u, u_{(1)}, u_{(2)}, \dots $, $V = W \bigcap Sym$.
In \cite{s} it is proved that for such a $V$ $v^{(q)} \leq v^{(1)} + q - 1$
for $q=1,2, \dots$ and $v^{(1)} \leq d+3$.

It is clear that for any ${\rm Q} \in V$ $\partial {\rm Q} / \partial y \in
V$.
In virtue of the Theorem \ref{t4} and the above formula (\ref{ansatzzz})
without loss of generality we may suppose that the characterystics of
the generalized symmetries of (\ref{eveq}) from $W$
of order not higher than $q$ has the form
\begin{equation} \label{evsym}
\eta = \exp(\lambda y) \sum_{j=0}^{v^{(q)}-1} \eta _{j}(u,u_{(1)}, \dots,
u_{(q)}) y^{j},
\end{equation}
where $\lambda \in \C$.

Since \cite{s} for any generalized symmetry ${\rm Q}$
of order $q \geq 2$ from $W$
\begin{equation} \label{qq}
\partial \eta /\partial u_{(q)} = {\rm const} (\partial G /\partial
u_{(d)})^{q/d}
\end{equation}
and $G$ is independant from $y$, the comparison of (\ref{evsym}) and
(\ref{qq})
shows that in such a case $\lambda=0$ in (\ref{evsym}), i.e. the generalized
symmetries of order $q \geq 2$ from $W$ depend on $y$ as polynomials of
order not higher than $v^{(q)} -1 \leq v^{(1)}+q-2$.

{\it Example 2.} Let us consider one-dimensional free Schr\"odinger
equation ($m=2$, $n=1$, $u_1 \equiv \psi$, $x_1 \equiv t, x_2 \equiv x$):
\begin{equation} \label{se}
L \psi \equiv (i \partial /\partial t + \partial ^{2}/\partial x^{2}
) \psi = 0 .
\end{equation}
and its {\it symmetry operators} \cite{fn} of the form
($\{ A, B \} \equiv A B +B A$, $p \equiv -i \partial /\partial x$)
\begin{equation} \label{so}
R = \sum _{j=0} ^{q} \{ \dots \{ h_j(x,t), \underbrace{p \},\dots
p \} }_{j~times} \equiv
\sum _{j=0} ^{q} a_j(x,t) \partial ^{j} /\partial x^{j},
\end{equation}
which should commute with $L$:
\begin{displaymath}
[R, L] \equiv R L - L R =0,
\end{displaymath}
what yields the following equations for $h_j$ \cite{bdn}:
\begin{equation} \label{deteq}
\begin{array}{lll}
\dot{h} _{j} = h'_{j-1},\quad j=1, \dots, q, \\
\dot{h} _{0} =0, \\ h'_{q} =0,
\end{array}
\end{equation}
where dot and prime denote partial derivatives with respect to $t$ and $x$.
The symmetry operators of the form (\ref{so}) of (\ref{se}) may be
considered as the elements of $Sym$ \cite{o}.

Let $W$ be the set of all the linear differential operators of the
form (\ref{so}).
It is known \cite{bdn} that for such a $W$ the dimensions $v^{(q)}$
of the subspaces $V^{(q)} = W \bigcap Sym^{(q)}$ (as usual,
$V=W \bigcap Sym$) are $v^{(q)} = (q+1)(q+2)/2$. Since
(\ref{se}) admits Lie symmetries $\partial
/\partial t$ and $\partial /\partial x$ and if $R$ is
symmetry operator of (\ref{se}), then so do
 $\partial R/\partial t$ and $\partial R/\partial x$, all the conditions
of the Theorem \ref{t4} are fulfilled. Using the formulas (\ref{ansatzzz}),
we obtain the following expression for the generic symmetry operator of
(\ref{se}) of order $q$:
\begin{equation} \label{ansatz}
R =\exp (\lambda t +\mu x) \sum_{k=0}^{v^{(q)}-1} \sum _{l=0}^{v^{(q)}-1}
\sum_{s=0}^{q} C_{kls} t^{k} x^{l} \partial ^{s} /\partial x^{s}.
\end{equation}
where $\lambda, \mu, C_{kls} \in \C$.
Hence, $h_{j}$ for $R$ (\ref{ansatz}) have the form:
\begin{equation} \label{hh}
h_{j} = \exp (\lambda t +\mu x) \sum_{k=0}^{v^{(q)}-1} \sum _{l=0}^{v^{(q)}-1}
 h_{jkl} t^{k} x^{l}, \quad h_{jkl} \in \C.
\end{equation}
The substitution of (\ref{hh}) into (\ref{deteq}) yields $\lambda
=\mu =0$ for $R \not \equiv 0$ (i.e. $h_j$ are polynomials
with respect to $t$ and $x$) and recurrent relations for $h_{jkl}$.
Thus, we partially recovered (by different means) the results of
\cite{bdn}. The results concerning the general solution of (\ref{deteq})
may be found there. Let us only
mention that $h_{j}$ is polynomial of order $j$ with
respect to $t$ and of order $q-j$ with respect to $x$.

The above examples show that our method is rather efficient tool for
the reduction of the equation (\ref{3}) for the coefficients of generalized
symmetries with respect to some "selected" variables,
which do not enter explicitly in the system of PDEs under consideration.

The requirement that $\partial {\rm Q}/\partial z_{A_{i}} \in V$,
$i=1, \dots, g$ for
${\rm Q} \in V$ in fact may be replaced by the weaker one:
there exist $g$ Lie vector fields ${\rm K}_{s} = \sum _{A =1} ^{m+n}
\omega _{A}^{(s)} (z) \partial /\partial z_{A}$, $s=1, \dots, g$,
such that $[{\rm K}_{i}, {\rm K}_{j}] =0$ for all $i,j=1, \dots, g$,
$[{\rm K}_{i}, {\rm Q}] \in V$ for any ${\rm Q} \in V$, $i=1, \dots, g$,
and in the generic point of $M^{(0)}$
${\rm rank } \parallel \omega _{A}^{(s)} \parallel _{A= \overline{1,
m+n}, s = \overline{1, g}} = g$. Really, in this case there exists
\cite{eis} such replacement of coordinates on $M^{(0)}$ (which, however, may
be defined only locally in each chart of $M^{(0)}$ but not globally)
$z \rightarrow z'$ that in new coordinates ${\rm K}_{s} = \partial /
\partial z'_{A_{s}}$, $A_{s} \in \lbrace 1, \dots, m+n \rbrace$, $s =1,
\dots, g$ and $[{\rm K}_{i}, {\rm Q}] =
\partial {\rm Q}/\partial z'_{A_{i}} \in V, i=1, \dots, g$ for any
${\rm Q} \in V$, i.e. we come back
to the situation, described in the Theorem \ref{t4}.

\section{Conclusions and Discussion}
The general idea of this work consists in the study of adjoint action of Lie
symmetries on $Sym$
(since $\partial {\rm Q} / \partial z_{A}
= [\partial _{A}, {\rm Q}]$ !).
The peculiarity of Lie symmetries is that the Lie bracket of the Lie symmetry
${\rm L}$ and some generalized symmetry ${\rm Q}$ of (\ref{1}) of order $q$ is
generalized symmetry of order not higher than $q$, i.e.
${\rm ad } _{L} \equiv [{\rm L}, \cdot] :
Sym ^{(q)} \rightarrow Sym ^{(q)}$ for any $q$. If the space $Sym ^{(q)}$
(or, in more general situation, $V^{(q)}$ ) is finite-dimensional for some
$q=q_1$, the operator ${\rm ad } _{L}$ posseses matrix representation
on it. This observation (for ${\rm L} = \partial _{A_{s}}, s=1,
\dots, g$) explains the analogy between the formulas (\ref {8aaa})
and the form of the ansatz for eigenfunctions of Hamiltonian in
quasiexactly solvable models \footnote{as pointed out to the author
Prof. R.Z. Zhdanov.}, since these
eigenfunctions belong to the common invariant space of the representation
of some Lie algebra by the first order linear differential operators, some
of which are just of the form $\partial /\partial z_{A}$ \cite{u}.

The above analogy with the theory of quasiexact solvability shows that
it would be very interesting to generalize the results of this paper to
the case of non-commutative algebra of Lie symmetries of (\ref{1}).
We intend to do it in further publications.

 I am sincerely grateful to Profs. A.G. Nikitin and
R.Z. Zhdanov for the fruitful discussion of the results of this work.

\end{document}